\documentclass[12pt]{amsart}
\usepackage{url}
\usepackage{epsfig}
\usepackage[all]{xy}
\usepackage{graphicx}
\usepackage{mathrsfs}

\newcommand{\defi}[1]{\textsf{#1}} 

\newcommand{\bA}{{\mathbb A}}
\newcommand{\bC}{{\mathbb C}}
\newcommand{\bF}{{\mathbb F}}
\newcommand{\bG}{{\mathbb G}}

\newcommand{\bP}{{\mathbb P}}
\newcommand{\bQ}{{\mathbb Q}}
\newcommand{\bR}{{\mathbb R}}
\newcommand{\bZ}{{\mathbb Z}}

\newcommand{\calA}{{\mathcal A}}

\newcommand{\calC}{{\mathcal C}}

\newcommand{\calE}{{\mathcal E}}

\newcommand{\calH}{{\mathcal H}}

\newcommand{\calL}{{\mathcal L}}

\newcommand{\calO}{{\mathcal O}}
\newcommand{\calP}{{\mathcal P}}
\newcommand{\calQ}{{\mathcal Q}}

\newcommand{\calU}{{\mathcal U}}
\newcommand{\calV}{{\mathcal V}}
\newcommand{\calW}{{\mathcal W}}
\newcommand{\calX}{{\mathcal X}}
\newcommand{\calY}{{\mathcal Y}}
\newcommand{\calZ}{{\mathcal Z}}

\newcommand{\wY}{\widetilde Y}

\DeclareMathOperator{\Tr}{Tr}

\DeclareMathOperator{\rk}{rk}
\DeclareMathOperator{\Char}{char}
\DeclareMathOperator{\inv}{inv}

\DeclareMathOperator{\Gal}{Gal}

\DeclareMathOperator{\Br}{Br}

\DeclareMathOperator{\Sym}{Sym}

\DeclareMathOperator{\Pic}{Pic}

\DeclareMathOperator{\Spec}{Spec}

\DeclareMathOperator{\Proj}{Proj}

\newcommand{\et}{{\rm \mathaccent 19 et}}

\newcommand{\ra}{\rightarrow}

\newtheorem{theorem}{Theorem}[section]
\newtheorem{lemma}[theorem]{Lemma}

\newtheorem{proposition}[theorem]{Proposition}

\theoremstyle{definition}

\newtheorem{example}[theorem]{Example}

\theoremstyle{remark}
\newtheorem{remark}[theorem]{Remark}
\newtheorem{remarks}[theorem]{Remarks}

\begin{document}

\title[Transcendental obstructions]{Transcendental obstructions to weak approximation on general K3 surfaces}
\author{Brendan Hassett}
\thanks{This research was supported by National Science Foundation Grants 0554491 and 0901645.}
\address{Department of Mathematics, Rice University, Houston, TX 77005, USA}
\email{hassett@rice.edu}
\urladdr{http://www.math.rice.edu/\~{}hassett}

\author{Anthony V\'arilly-Alvarado}
\address{Department of Mathematics, Rice University, Houston, TX 77005, USA}
\email{varilly@rice.edu}
\urladdr{http://www.math.rice.edu/\~{}av15}

\author{Patrick Varilly}
\address{Department of Chemistry, University of California at Berkeley, Berkeley, CA 94720}
\email{patvarilly@gmail.com}
\urladdr{http://gold.cchem.berkeley.edu/\~{}pvarilly}

\begin{abstract}
We construct an explicit K3 surface over the field of rational numbers that has geometric Picard rank one, and for which there is a transcendental Brauer-Manin obstruction to weak approximation. To do so, we exploit the relationship between polarized K3 surfaces endowed with particular kinds of Brauer classes and cubic fourfolds.
\end{abstract}

\maketitle

\section{Introduction}\label{S:introduction}

Let $X$ be a smooth projective geometrically integral variety over a number field $k$. Assume that the set $X(\bA_k)$ of adelic points is nonempty; we say that $X$ satisfies \defi{weak approximation} if the natural embedding $X(k) \hookrightarrow X(\bA_k)$ is dense for the adelic topology.  Following Manin, we may use the Brauer group $\Br(X) := H^2_\et(X,\bG_m)$ to construct a set $X(\bA_k)^{\Br} \subseteq X(\bA_k)$ such that the closure of $X(k)$ in $X(\bA_k)$ already lies in $X(\bA_k)^{\Br}$. If $X(\bA_k) \setminus X(\bA_k)^{\Br}\neq \emptyset$, then we say there is a \defi{Brauer-Manin obstruction to weak approximation} on $X$. In this way, the group $\Br X$ may be used to obtain qualitative information about the distribution of rational points on $X$.

Fix an algebraic closure $\bar k$ of $k$, and let $\overline{X} = X\times_k\bar k$. Write $\Br_0(X)$ for the image of $\Br k$ in $\Br (X)$, and let $\Br_1 (X) = \ker \big(\Br (X) \to \Br(\overline{X})\big)$. An element of $\Br (X)$ is \defi{algebraic} if it lies in $\Br_1 (X)$, \defi{transcendental} otherwise. 

The first examples of transcendental Brauer elements in an arithmetic context are due to Harari~\cite{Harari}.  Rational surfaces posses no such elements because the Brauer group is a birational invariant, and $\Br \big(\bP^2_{\bar k}\big) = 0$. In contrast, for K3 surfaces the group $\Br (X)/\Br_1 (X)$ can be nontrivial, albeit finite~\cite[Theorem 1.2]{SkorobogatovZarhin}. For these surfaces, the arithmetically interesting groups $\Br (X)/\Br_0 (X)$, $\Br (X)/\Br_1 (X)$ and $\Br_1 (X)/\Br_0 (X)$ are the object of much recent research~\cite{Bright, SkorobogatovZarhin, IeronymouSkorobogatovZarhin, KreschTschinkel, LoganVanLuijk, SkorobogatovZarhin2, Corn}. Explicit transcendental elements of $\Br (X)$, or the lack thereof, play a central role in \cite{Wittenberg,SkorobogatovSwinnertonDyer,HarariSkorobogatov,Ieronymou}. In all cases, the K3 surfaces considered are endowed with an elliptic fibration, and this additional structure is essential to the computation of transcendental classes. A `general' K3 surface (i.e., one with geometric Picard rank one), however, does not carry such a structure. With this in mind, our main result is as follows.

\begin{theorem}
\label{thm: Main}
Let $X$ be the K3 surface of degree $2$ given by
\begin{equation}
w^2 = \det
\begin{pmatrix}
2(2x + 3y + z) & 3x + 3y & 3x + 4y & 3y^2 + 2z^2 \\
3x + 3y & 2z & 3z & 4y^2 \\
3x + 4y & 3z & 2(x+3z) & 4x^2 + 5xy + 5y^2 \\
3y^2 + 2z^2 & 4y^2 & 4x^2 + 5xy + 5y^2 & 2(2x^3 + 3x^2z + 3xz^2 + 3z^3)
\end{pmatrix}
\label{eq:theK3}
\end{equation}
in $\bP_\bQ(1,1,1,3)$. Then $X$ has geometric Picard rank one, and there is a transcendental Brauer-Manin obstruction to weak approximation on $X$. The obstruction arises from a quaternion Azumaya algebra $\mathscr{A} \in \Br \big(k(X)\big)$ in the image of the natural injection $\Br (X) \hookrightarrow \Br \big(k(X)\big)$. Explicitly, we have $\mathscr{A} = (\alpha,\beta)$, where
\begin{align*}
\alpha &= \frac{9x^2 + 18xy - 8xz + 9y^2 - 12yz - 4z^2}{4(2x + 3y + z)^2},\ \text{and} \\
\beta &= -\frac{9x^3 + 18x^2y + x^2z + 9xy^2 + 3xyz - 10xz^2 + 7y^2z - 9yz^2 - 3z^3}{(2x + 3y + z)(9x^2 + 18xy - 8xz + 9y^2 - 12yz - 4z^2)}.
\end{align*}
\end{theorem}

The pair $(X,\mathscr{A})$ in Theorem~\ref{thm: Main} is naturally associated with the smooth cubic fourfold $Y \subseteq \bP^5_\bQ$ given by
\begin{equation}
\label{eq:fourfold}
\begin{split}
2X_1^2Y_1 &+ 3X_1^2Y_2 + X_1^2Y_3 + 3X_1X_2Y_1 + 3X_1X_2Y_2 + 3X_1X_3Y_1 + 4X_1X_3Y_2 + 3X_1Y_2^2 \\
&+ 2X_1Y_3^2 + X_2^2Y_3 + 3X_2X_3Y_3 + 4X_2Y_2^2 + X_3^2Y_1 + 3X_3^2Y_3 + 4X_3Y_1^2 + 5X_3Y_1Y_2 \\
&+ 5X_3Y_2^2 + 2Y_1^3 + 3Y_1^2Y_3 + 3Y_1Y_3^2 + 3Y_3^3 = 0.
\end{split}
\end{equation}
This fourfold contains the plane $\{Y_1 = Y_2 = Y_3 = 0\}$.  Quite generally, we show that, for any field $k$, we may construct, from a given smooth cubic fourfold $Y \subseteq \bP^5_k$ containing a plane $P$ and no other plane $P'$ meeting $P$ along a line, a K3 surface of degree $2$ over $k$, together with a $2$-torsion element of $\Br(X)$; see Theorem~\ref{theorem:getAzumaya}.  This construction is well-known in the case $k = \bC$~\cite{Voisin,vanGeemen}; Voisin used it in her
proof of the Torelli theorem for cubic fourfolds.  The main result of \cite{Has99} may be interpreted as saying that $Y$ is rational when the $2$-torsion element of $\Br(X)$ is trivial.  Recent work of Macri and Stellari \cite{MS} recasts these constructions in the context of derived categories.

There is no \emph{a priori} guarantee that the K3 surface furnished by the construction above has geometric Picard rank one.  We use a method of van Luijk, together with the work of Elsenhans and Jahnel on K3 surfaces of degree $2$, to assemble a K3 surface for which we can prove that this is the case~\cite{vanLuijk,ElsenhansJahnelrankone,ElsenhansJahneldoublecover}.
In essence, we find smooth cubic fourfolds $\mathfrak{C}_2$ and $\mathfrak{C}_3$ over the finite fields $\bF_2$ and $\bF_3$, respectively, whose corresponding K3 surfaces have rank-two geometric Picard lattices, with discriminants in different square classes (i.e., whose quotient is not a square). Then the K3 surface arising from a cubic fourfold whose reductions at $2$ and $3$ coincide with $\mathfrak{C}_2$ and $\mathfrak{C}_3$, respectively, has geometric Picard rank one. This explains why the coefficients in~\eqref{eq:fourfold}, as well as the entries of the matrix in~\eqref{eq:theK3}, are all integers less than six.


\subsection{Outline of the Paper} In \S\ref{S:hodgetheory}, we explain the known relation between special cubic fourfolds over $\bC$ and complex projective polarized K3 surfaces furnished with particular kinds of $2$-torsion Brauer classes, using Hodge theory. In \S\ref{s:quadrics}, we generalize some well-known facts about quadratic forms over base fields to the case of arbitrary base schemes. 
In \S\ref{S:CubicToK3}, we spell out how to explicitly construct a K3 surface from a special cubic fourfold over a field $k$, paying particular attention to the case where $\Char(k) = 2$. 
We apply the results of \S\ref{s:quadrics} in \S\ref{S:constructing} to give geometric proofs, valid over any field, of the constructions outlined in \S\ref{S:hodgetheory}.
We use these constructions over fields of arithmetic interest in the rest of the paper. 
In \S\ref{S:Picard}, we use van Luijk's method to construct a K3 surface $X$ over $\bQ$, arising from a special cubic fourfold, with geometric Picard rank one. In \S\ref{S: proof}, we compute the unramified quaternion algebra in $\Br(\bQ(X))$ associated to our K3 surface by the constructions of \S\ref{S:constructing}, and use it to prove Theorem~\ref{thm: Main}. Finally, in \S\ref{Section: computations}, we describe how we implemented our construction on a computer, and we outline a few speed-ups that made the project possible.

\subsection{Notation}
Throughout $k$ denotes a field (precise hypotheses are given when needed), and $\bar k$ is a fixed algebraic closure of $k$. If $X$ and $Y$ are $S$-schemes, then $X_Y := X\times_S Y$. If $Y = \Spec R$, then we write $X_R$ instead of $X_{\Spec R}$; if $X$ is a $k$-scheme, we write $\overline{X}$ for $X_{\bar k}$. For an integral scheme $X$ over a field, we write $k(X)$ for the function field of $X$. When $X$ is a K3 surface over a field, we write $\rho(X)$ for the rank of the free abelian group $\Pic(\overline{X})$.
 
Given a Hodge structure $H$ of weight $n$, we write $H(-1)$ for its Tate twist,
the Hodge structure with an identical filtration but of weight $n+2$.  The Grassmannian of
$r$-dimensional subspaces of $\bP^n$ is denoted $\bG(r,n)$.  

\subsection*{Acknowledgements}  
\begin{figure}[h]
\begin{center}
\includegraphics[scale=0.5]{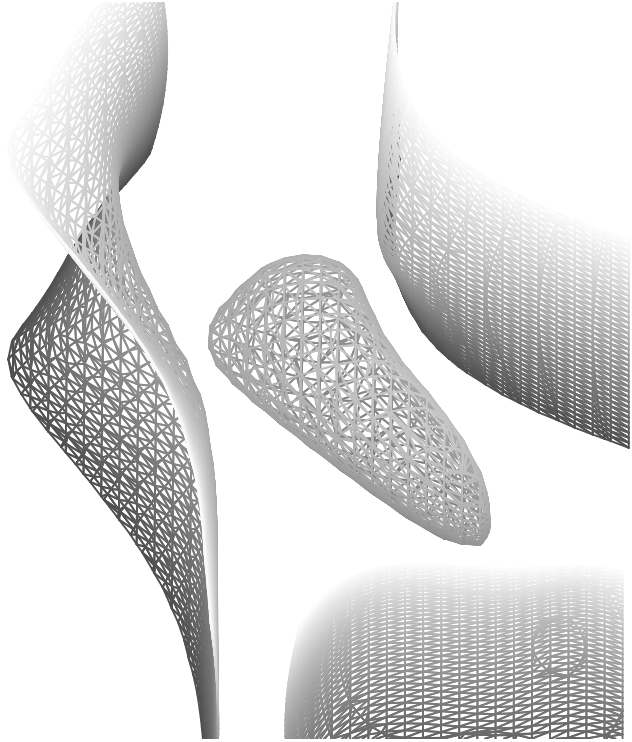}
\end{center}
\caption{The real component (center) of the K3 surface of Theorem~\ref{thm: Main} where the real invariant of $\calA$ is nontrivial (view of the $z = 1$ affine patch)}
\label{figure:real points}
\end{figure}

We are grateful to Andreas-Stephan Elsenhans, David Harari, Kelly McKinnie, Yuri Tschinkel and Olivier Wittenberg for discussions about this work.  
Many computations presented here were done with \texttt{Magma}~\cite{BCP}.  Figure~\ref{figure:real points} was created using the Maple software package.

\section{Hodge theory and two-torsion elements of the Brauer group}
\label{S:hodgetheory}
Given a smooth cubic fourfold $Y\subseteq \bP^5_\bC$, containing a plane $P$ that does not meet any other plane, we may construct a K3 surface $X$ as follows.  Let $\pi \colon\widetilde{Y} \to \bP^2_\bC$ be the quadric bundle induced by projection away from $P$. Then there is a smooth sextic curve $C \subseteq \bP^2_\bC$ that parametrizes the singular fibers of $\pi$~\cite[Lemme~2, p.\ 586]{Voisin}. The K3 surface $X$ is the double cover of $\bP^2_\bC$ ramified along $C$.

In this section, we explain how to go in the opposite direction in the presence of extra data. To wit, following~\cite{Voisin,vanGeemen}, we can recover a smooth cubic fourfold $Y\subseteq \bP^5_\bC$ (containing a plane) from a pair $\big((X,f),\alpha\big)$, consisting of a polarized complex projective K3 surface $(X,f)$ such that $\Pic X = \bZ f$, and a particular kind of element $\alpha \in \Br(X)[2]$ (the $2$-torsion subgroup of the Brauer group). This construction is Hodge-theoretic. 

This point of view has arithmetic implications:  Given a two-torsion element of the Brauer group of a degree-two K3 surface, all defined over a number
field, it is natural to seek a canonically-defined unramified Azumaya algebra, arising from a concrete geometric construction, representing the Brauer class. 
We produce explicit examples of this type, and apply the resulting Brauer element to weak approximation questions.
While the content of this section is not logically necessary for our results, it illustrates the geometric motivations behind our 
approach. \\

Let $X$ be a complex projective K3 surface, and let $\left<\,,\right>$ denote
the cup product pairing on its middle cohomology.  Consider the cohomology 
group $H^2_{\et}(X,\mu_2)$, where $\mu_2=\{\pm 1 \}$.  
We have in mind the \'etale topology, but standard comparison theorems 
\cite[III.3.12]{Milne} allow us
to identify this with $H^2(X(\bC),\mu_2)\simeq H^2(X(\bC),\bZ/2\bZ)$, 
sheaf cohomology for the underlying analytic space.
In particular, we have 
$$H^2_{\et}(X,\mu_2) \simeq (\bZ/2\bZ)^{22}.$$
The cup product map
$$
\begin{array}{rcl}
  H^2_{\et}(X,\mu_2) \times H^2_{\et}(X,\mu_2) & \ra &  H^4_{\et}(X,\mu_2^{\otimes 2}) \simeq \bZ/2\bZ \\
        (\alpha,\beta) & \mapsto & \left<\alpha,  \beta\right>
\end{array}
$$
is nondegenerate.  Since the intersection form on the middle
cohomology of a K3 surface is even, we have
$\left<\alpha,\alpha\right>=0$
for each $\alpha \in H^2_{\et}(X,\mu_2)$.  Nevertheless, there
exists a quadratic form 
$$q_{\left<,\right>}\colon H^2_{\et}(X,\mu_2) \ra H^4_{\et}(X,\mu_2^{\otimes 2})$$
constructed as follows:  Take the integral (or $2$-adic) intersection quadratic
form on the middle cohomology, divide this by two, and reduce modulo two.  
An element $\alpha \in H^2_{\et}(X,\mu_2)$ is {\em even} or {\em odd}
depending on the parity of $q_{\left<,\right>}(\alpha)$.

Consider the exact sequence
$$\Pic(X) \stackrel{\times 2}{\ra}
\Pic(X) \ra H^2_{\et}(X,\mu_2) \ra \Br(X)[2] \ra 0,$$
where the first arrow is the mod $2$ cycle-class map.
If $\rho(X)=\rk(\Pic(X))$ then 
$$\Br(X)[2] \simeq (\bZ/2\bZ)^{22-\rho(X)}.$$

Suppose that $X$ has degree $2$ and is endowed with a polarization $f$.
In particular, $f$ is an ample divisor with $\left<f,f\right>=2$.  
Consider the subset
$$\{\alpha \in H^2_{\et}(X,\mu_2): q_{\left<,\right>}(\alpha)\equiv \left<\alpha,f\right>\equiv 1\pmod{2} \},$$
i.e., the odd non-primitive classes.  Note that $\alpha$ is in this set if and only if $\alpha+f\pmod{2}$
is as well, so we can consider
$$
B(X,f):=\{\alpha \in H^2_{\et}(X,\mu_2)/\left<f\right> : q_{\left<,\right>}(\alpha)\equiv \left<\alpha,f\right>\equiv 1\pmod{2} \}.$$
The map $B(X,f) \ra \Br(X)[2]$ is injective  if and only if $\Pic(X)=\bZ f$.

Voisin \cite[\S 1]{Voisin} and van Geemen \cite[\S 9.7]{vanGeemen} offer a geometric construction
for Severi-Brauer varieties smooth of relative dimension one over $X$, indexed by $B(X,f)$.  
We summarize their results as follows.
\begin{proposition} \label{prop:vGV}
Let $(X,f)$ be a complex projective K3 surface of degree two with $\Pic(X)=\bZ f$ and let $\alpha \in B(X,f)$.  
Then there exists a unique cubic fourfold $Y$ containing a plane $P$, with variety of lines $F$, and a correspondence
\begin{equation*}
\xymatrix{
 & \calW \ar[dl]_{\pi_1}\ar[dr]^{\pi_2} & & \calZ \ar[dl]_{p}\ar[dr]^{q} & \\
X & & F & & Y
}
\end{equation*}
where $\calW \subset F$ parametrizes the lines in $Y$ incident to $P$,
$\pi_1$ is a smooth $\bP^1$-bundle (in the \'etale topology), $\pi_2$ is an inclusion, and
$$\calZ=\{[\ell,y]: y\in \ell \subset Y \}$$ 
is the incidence correspondence between $Y$ and $F$.   Moreover, $[\calW]=\alpha$ in the Brauer
group of $X$.  Precisely, the correspondence induces an isogeny of Hodge structures
$$\begin{array}{ccc}
H^4(Y,\bZ)  & & H^2(X,\bZ)(-1) \\
\cup & & \cup \\
j\colon\{h^2,P\}^{\perp} & \hookrightarrow &  f^{\perp},
\end{array}
$$
where $h$ is the hyperplane class on $Y$, whose image  equals 
\[
\{\eta \in H^2(X,\bZ): \left<\eta,f\right>=0, \left<\eta,\alpha\right>\equiv 0 \pmod{2}\}.
\]  
The cubic fourfold $Y$ does not contain another plane $P'$ that meets $P$. 
\end{proposition}
We sketch the construction of $j$:  The incidence correspondence $\calZ$ induces an Abel-Jacobi map 
\[
p_*q^*\colon H^4(Y,\bZ) \ra H^2(F,\bZ)(-1),
\]
i.e., an isogeny of Hodge structures 
that is an isomorphism between the primitive cohomology groups \cite{BD}.  In particular, we may regard
\[
\{h^2,P\}^{\perp}  \subset H^2(F,\bZ)(-1).
\]
The pullback maps
\[
\pi_1^*\colon H^2(X,\bZ) \ra H^2(\calW,\bZ),  \quad \pi_2^*\colon H^2(F,\bZ) \ra H^2(\calW,\bZ),
\]
are injective;  indeed, for the former we know that $\pi_1$ is projective and for the latter we use
\cite[Proposition 1, p.\ 582]{Voisin}.  Thus we can compare 
\[
\pi_1^* \big( f^{\perp}\big)\quad\text{and}\quad \pi_2^*p_*q^* \left(\{h^2,P\}^{\perp}\right)(1)
\]
inside $H^2(W,\bZ)$, and find that the latter is an index-two subgroup of the former \cite[Lemme 5, p.\ 583]{Voisin}.

\begin{remarks}\ 
\begin{enumerate}
\item
The geometric manifestation of $\alpha$ by $\calW$ in~Proposition~\ref{prop:vGV} is the kind of presentation for Brauer elements that we emulate in arithmetic contexts; see~\S\ref{S:constructing}.
\item
The hypothesis that $\Pic(X)=\bZ f$ in Proposition~\ref{prop:vGV} is much
stronger than is necessary.
\item
The analysis above does not apply to all $\alpha \in \Br(X)[2]$;  van Geemen \cite[\S 9]{vanGeemen}
offers geometric interpretations for the classes we have not considered.
\item
Macri and Stellari have recast this discussion in the language of $\alpha$-twisted sheaves on $X$~\cite{MS}.
\end{enumerate}
\end{remarks}

\section{Quadratic forms and maximal isotropic subspaces}
\label{s:quadrics}
To generalize the constructions of \S\ref{S:hodgetheory}, we first generalize a few well-known results on
quadratic forms over fields to the case of arbitrary base schemes. We include proofs for lack of a suitable
reference.

\subsection{Quadratic forms over fields}
We begin with some elementary geometric observations,
which can be found in \cite[Chapter~22]{Harris}, at least over
fields of characteristic different from two.  A characteristic-independent
approach can be found in \cite[\S\S 7-14, 85]{EKM}. Fix $n\ge 2$ and let $Q \subset \bP_k^{2n-1}$ and $Q'\subset \bP_k^{2n-2}$
denote smooth quadric hypersurfaces over an algebraically closed field $k$.  Then $Q$
is isomorphic to 
$$x_1x_2+x_3x_4+\ldots+x_{2n-1}x_{2n}=0,$$
and $Q'$ is isomorphic to
$$x_1x_2+x_3x_4+\ldots+x_{2n-1}^2=0$$
\cite[\S 7]{EKM}.  Consider the schemes
$$W \subset \bG(n-1,2n-1), \quad W'\subset \bG(n-2,2n-2)$$ 
of maximal isotropic subspaces of $Q$ and $Q'$, respectively.  
Both schemes are smooth projective of dimension 
$$n^2-\binom{n+1}{2}=(n-1)n-\binom{n}{2}=\binom{n}{2};$$
furthermore, $W'$ is irreducible and $W=W_1 \cup W_2$, where
each component 
$W_i\simeq W'$
\cite[\S 85]{EKM}.
Indeed, if we fix a realization of $Q'$ as a hyperplane section of $Q$
$$Q' \hookrightarrow Q$$
then we have a restriction morphism
$$\begin{array}{rcl}
\rho\colon W & \ra & W' \\
\Lambda & \mapsto & \Lambda \cap Q'
\end{array}
$$
mapping each component of $W$ isomorphically onto $W'$.  For example,
if 
$$\Lambda=\{x_1=x_3=\cdots=x_{2n-1}=0 \} \subset Q$$
then 
$$\rho(\Lambda)=\{x_1=x_3=\cdots=x_{2n-1}=0 \} \subset Q'$$
and 
$$\rho^{-1}(\rho(\Lambda))=\left\{\Lambda, \bar{\Lambda}=\{x_1=\cdots=x_{2n-3}=x_{2n}=0 \}\right\}.$$

Suppose now that $Q\subset \bP_k^{2n-1}$ is a quadric over an arbitrary field $k$.  
Consider the \defi{Clifford algebra} $C(Q)$ over $k$ \cite[\S 11]{EKM}.
It decomposes into even and odd parts:
\[
C(Q)=C_0(Q) \oplus C_1(Q).
\]
When $Q$ is smooth, we have:
\begin{itemize}
\item{$C(Q)$ is central simple;}
\item{$C_0(Q)$ is separable, and its center is a quadratic \'etale
$k$-algebra, which is denoted by $Z(Q)$ and is called the \defi{discriminant algebra} \cite[\S 13]{EKM}}
\item{$C_0(Q)$ is central simple as an algebra over $Z(Q)$.}
\end{itemize}
When $\Char(k)\neq 2$, we can diagonalize
$$Q=\{a_1x_1^2+\cdots+a_{2n}x_{2n}^2 =0 \}$$
and $Z(Q)=k[w]/\left<w^2-c\right>$ with $c=(-1)^n a_1 \ldots a_{2n}.$
When $\Char(k)=2$, we can write 
$$Q=\{a_1x_1^2+x_1x_2+b_1x_2^2+\cdots+ a_n x_{2n-1}^2+x_{2n-1}x_{2n}+b_nx_{2n}^2=0 \}$$
and $Z(Q)=k[w]/\left<w^2+w+c\right>$ where $c=a_1b_1+\cdots+a_nb_n$, the \defi{Arf invariant}.  
Finally, if $W$ denotes the variety of maximal isotropic subspaces of $Q$, then
we have \cite[\S 85]{EKM}:
\begin{itemize}
\item{the geometric components of $W$ are defined over $Z(Q)$;}
\item{when $n=2$, each component of $W$ over $Z(Q)$ is the conic associated
with the quaternion algebra $C_0(Q)$ over $Z(Q)$.}
\end{itemize}

Assume again that $k$ is algebraically closed, and 
suppose that $Q_0 \subset \bP_k^{2n-1}$ is a quadric
with an isolated singularity.  In suitable coordinates, it has equation
$$x_1x_2+x_3x_4+\cdots+x_{2n-1}^2=0.$$
i.e., it is a cone over $Q'$.  The scheme
$W_0 \subset \bG(n-1,2n-1)$ of $(n-1)$-dimensional isotropic subspaces
is thus $W'$, set-theoretically.  Note that $W_0$ is Cohen-Macaulay and
thus has no embedded points---indeed, it is
the locus where a section of a rank-$\binom{n+1}{2}$ vector
bundle on $\bG(n-1,2n-1)$ vanishes.  However, $W_0$ has multiplicity
two at its generic point.  

\subsection{Quadratic forms over schemes}

We shall generalize the above facts to arbitrary base schemes.  Let $\calH$ be the Hilbert scheme parametrizing quadratic forms in $2n$ variables
$$\sum_{1\le i \le j \le 2n} a_{ij}X_iX_j.$$
Concretely, we have
$$\calH=\bP^{\binom{2n+1}{2}}_{\bZ}=\Proj \bZ[a_{11},\ldots,a_{2n\, 2n}].$$  
Consider the open subsets 
$$\calV \subset \calU \subset \calH,$$
corresponding to quadrics
that are geometrically smooth and have at most one isolated singularity.  
In characteristics $\neq 2$, $\calU$ (resp.\ $\calV$) corresponds to quadrics whose associated
symmetric matrix has rank $\ge 2n-1$ (resp.\ $2n$).  

Let $\Delta \subset \calH$ denote the \defi{discriminant divisor}, defined as follows:  
If $q\colon\calQ \ra \calH$ is the universal family of quadrics then the 
differential of $q$ drops rank along a subscheme $\Sigma \subset \calQ$;  
$\Delta$ is the divisor obtained by pushing forward $\Sigma$ via $q$.  
Note that the discriminant {\em algebra} ramifies over the discriminant
{\em divisor}.  
\begin{example} \label{example:nosection}
In characteristic $2$, the morphism from $\Sigma$ to its set-theoretic
image is purely inseparable of degree two.  
In particular, $\Delta$ has multiplicity $2$ in characteristic $2$.  
Indeed, given the family
$$x_1x_2+x_3^2+sx_4^2=0$$
of non-smooth quadrics over the $s$-line, the singular locus is 
$$\{x_1=x_2=x_3^2+sx_4^2=0\}$$
which is purely inseparable over the $s$-line.  
There is no section through the singularities.  
\end{example}

\begin{example}
\label{ex:discriminants}
When $n=1$ we have 
$$\Delta=\{-4a_{11}a_{22}+a_{12}^2=0 \};$$
the corresponding quadratic algebra is
\[
Z(Q)= \begin{cases}
k[w]/\left< w^2-a_{12}^2+4a_{11}a_{22}\right> & \Char(k) \neq 2, \\
k[w]/\left<w^2+wa_{12}+a_{11}a_{22}\right> & \Char(k)=2.
\end{cases}
\]
Note the first form specializes to the second
via the substitution
$w\mapsto 2w+a_{12}$ and division by four.  

For $n=2$, setting
$L=a_{12}a_{34}+a_{13}a_{24}+a_{23}a_{14}$, 
$N=a_{11}a_{22}a_{33}a_{44}$, and
\[
\begin{split}
M = &-(a_{12}a_{23}a_{34}a_{14} + a_{13}a_{23}a_{24}a_{14}+a_{12}a_{24}a_{13}a_{34}) \\
    &+(a_{11}a_{23}a_{24}a_{34}+a_{22}a_{13}a_{34}a_{14}+a_{33}a_{12}a_{24}a_{14}+a_{44}a_{12}a_{23}a_{13}) \\
    &-(a_{11}a_{22}a_{34}^2+a_{11}a_{33}a_{24}^2+a_{11}a_{44}a_{23}^2+a_{22}a_{33}a_{14}^2+a_{22}a_{44}a_{13}^2+a_{33}a_{44}a_{12}^2),
\end{split}
\]
we have
$$\Delta=\{L^2+4M+16N=0\}.$$
Note that the even coefficents are divisible by four.  
The corresponding quadratic algebra is
\[
Z(Q)=\begin{cases}
k[w]/\left<w^2-(L^2+4M+16N)\right> & \Char(k)\neq 2, \\
k[w]/\left<w^2+wL+M\right> & \Char(k)=2.
\end{cases}
\]
The first form specializes to the second via the substitution
$w\mapsto 2w+L.$
\end{example}

Let $\calW \subset \bG(n-1,2n-1) \times \calU$ parametrize $(n-1)$-dimensional
isotropic subspaces in the projective hypersurfaces parametrized by $\calU$,
and let $r\colon\calW \ra \calU$ be the structure map.   Since $\calW$ is an open
subset or a projective bundle over $\bG(n-1,2n-1)$, it is smooth over $\Spec(\bZ)$.  

\begin{proposition} \label{proposition:Stein}
The Stein factorization 
$$r\colon\calW \stackrel{\pi}{\ra} \calX \stackrel{\phi}{\ra} \calU$$
is the composition of a smooth morphism by a flat double cover
branched along $\Delta$.  The geometric fibers of $\pi$ are all isomorphic to $W'$.
\end{proposition}
\begin{proof}
The fact that the Stein factorization is a double cover branched along $\Delta$
follows from our previous discussion, i.e., that the components of the maximal
isotropic subspaces are defined over the discriminant algebra, which ramifies
over the discriminant divisor.  The assertion on the geometric fibers of $\pi$
follows from the isomorphism $W_i\simeq W'$ above, once we know that $\pi$ is
smooth.  We therefore focus on this last statement.

We first show that $\phi$ is flat.
The scheme $\calW$ is defined as the locus where a section of a rank-$\binom{n+1}{2}$ vector bundle
on the Grassmannian vanishes, with the expected dimension.  Hence $\calW$ is flat over $\calU$.  
The geometric fibers of $r$ are homogeneous 
varieties whose structure sheaves have no higher cohomologies.  Thus 
$\bR^ir_*\calO_{\calW}=0$ for each $i>0$
and $r_*\calO_{\calW}$ is flat by cohomology-and-base-change.  

Consider the pull-back of the universal quadric to the discriminant cover/Stein
facorization
\[
q''\colon\calQ'':=\calQ \times_{\calU} \calX \ra \calX;
\]
$\calW$ tautologically determines a component of the variety
of maximal isotropic subspaces of $\calQ''$.  
We claim that there exists a smooth subscheme $\Sigma'' \subset \calQ''$ containing 
the singularities of the fibers and meeting each singular fiber in a subscheme
of length one.  In other words, over the reduced discriminant divisor in $\calX$ there
exists a section through the singularities (cf. Example~\ref{example:nosection}).

This statement is clear except in characteristic two, where the non-smooth
points may be defined over a purely inseparable extension.  
It suffices to exhibit the sections after \'etale localization on the base.  After such localization,
our families of quadrics have normal form
\[
x_1x_2+x_3x_4+\cdots+x_{2n-1}^2+tx_{2n-1}x_{2n}+sx_{2n}^2=0;
\]
indeed, this is versal for deformations of 
\[
x_1x_2+x_3x_4+\cdots+x_{2n-1}^2=0.
\]
Here, the discriminant algebra takes the form
\[
z^2+tz+s=0.
\]
More symmetrically, we write
\[
\alpha\beta=s \quad\text{and}\quad \alpha+\beta=t,
\]
so that $\calQ''$ is locally of the form
\[
x_1x_2+x_3x_4+\cdots+(x_{2n-1}+\alpha x_{2n})(x_{2n-1}+\beta x_{2n})=0
\]
with discriminant locus $\alpha=\beta$.  
Here we can write 
\[
\Sigma''=\{x_1=x_2=\cdots =x_{2n-2}=x_{2n-1}+\alpha x_{2n}=x_{2n-1}+\beta x_{2n}=0\}.
\]

The local analysis yields additional dividends:  We see that $\calQ''$ has codimension-$(2n-1)$ ordinary
double points along $\Sigma''$.  Consider the modification
$$\beta\colon\widetilde{\calQ}:=\mathrm{Bl}_{\Sigma''}\calQ'' \ra \calQ''$$
and the induced morphism
$$\tilde{q}\colon\widetilde{\calQ} \ra \calX,$$
which remains flat over $\calX$;  let $\calE$ denote its exceptional divisor.  Fix a point $p$ in the discriminant on $\calX$, so that
$\calQ''_p={q''}^{-1}(p)$ is singular at $\sigma=\Sigma''_p$.  Then we have
$$\widetilde{\calQ}_p:=\tilde{q}^{-1}(p)=Q \cup_{Q'} \mathrm{Bl}_{\sigma}\calQ''_p,$$
where $Q=\calE \cap \widetilde{\calQ}_p$ is a smooth quadric of dimension $2n-2$, $Q'\subset Q$ is a smooth hyperplane
section, and $Q'\subset \mathrm{Bl}_{\sigma}\calQ''_p$ is the exceptional divisor.  
Note that $\mathrm{Bl_{\sigma}}\calQ''_p$ is a projective bundle over a smooth quadric of 
dimension $2n-3$, with section $Q'$.  

\begin{figure}
\begin{center} 
\includegraphics{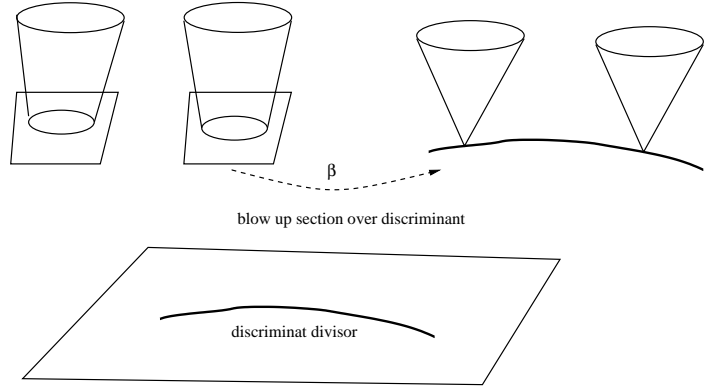}
\end{center}
\caption{First modification in the $n=2$ case}
\end{figure}

Regard $\calQ'' \subset \calX \times \bP^{2n-1}$ and let $\calO_{\calQ''}(1)$ denote the pull back of $\calO_{\bP^{2n-1}}(1)$ 
via projection.  Consider the invertible sheaf
$$\beta^*\calO_{\calQ''}(1)\otimes \calO_{\widetilde{\calQ}}(-\calE),$$
which is semiample and induces a contraction 
\[
\xymatrix{
\widetilde{\calQ} \ar[rd]_{\tilde{q}} \ar[rr]^\gamma & &  \widehat{\calQ} \ar[ld]^{\hat q}\\
 & \calX &
}
\]
This has the following properties:
\begin{itemize}
\item{$\gamma$ is the blow up along a smooth codimension-two subvariety $\calY\subset \widehat{\calQ}$,
supported over the discriminant;  $\calY$ intersects singular fibers of $\hat{q}$ in smooth
hyperplane sections;}
\item{for $p$ in the discriminant, $\gamma|\widetilde{\calQ}_p$ is 
an isomorphism on $Q$, and it collapses $\mathrm{Bl}_{\sigma}\calQ''_p$ along its ruling to $Q'\subset Q$;}
\item{$\hat{q}$ is a smooth quadric bundle;}
\item{the variety $\widehat{\calW} \ra \calX$ parametrizing maximal isotropic subspaces in fibers of $\hat{q}$ is
smooth over $\calX$ with two connected components.}
\end{itemize}

\begin{figure}
\begin{center} 
\includegraphics{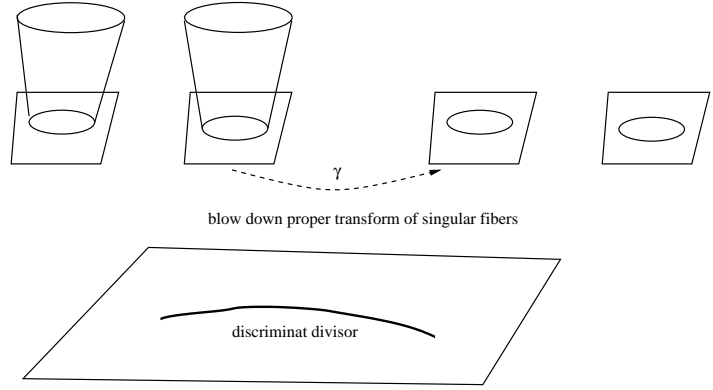}
\end{center}
\caption{Second modification in the $n=2$ case}
\end{figure}

\begin{example}
Given the quadric 
$$\calQ=\{x_1x_2+\cdots+x_{2n-1}^2 - \alpha^2 x_{2n}^2=0\}$$
over the complex affine $\alpha$-line, we have
$$\widehat{\calQ}=
\{x_1x_2+\cdots+x_{2n-1}^2 - y^2=0\},$$
where $y=\alpha x_{2n}$.
\end{example}

\begin{lemma}
Consider the part of the relative Hilbert scheme of $\tilde{q}\colon\widetilde{\calQ}\ra \calX$
$$\widetilde{\calW}=\widetilde{\calW}_1 \cup \widetilde{\calW}_2 \ra \calX$$
that, for smooth fibers of $\tilde{q}$, parametrizes maximal isotropic subpaces.  
These satisfy the following properties:
\begin{itemize}
\item{for $p$ in the discriminant, $\widetilde{\calW}_p$ parametrizes 
$$\{\Lambda:\Lambda \subset \widetilde{\calQ}_p=Q \cup_{Q'} \mathrm{Bl}_{\sigma}\calQ''_p\},$$
where $\Lambda \cap Q$ and $\Lambda \cap Q'$ are maximal isotropic, and 
$\Lambda \cap \mathrm{Bl}_{\sigma}\calQ''_p$ is the proper transform of a maximal isotropic subspace
in $\calQ''_p$;}
\item{$\gamma$ induces an isomorphism $\widetilde{\calW}\stackrel{\sim}{\ra} \widehat{\calW}$;}
\item{$\beta$ induces an isomorphism between one of the $\widetilde{\calW}_i$ 
and $\calW$.}
\end{itemize}
\end{lemma}
\begin{proof} (of the Lemma)
The morphism $\gamma$ blows up a smooth codimension-two subvariety $\calY \subset \widehat{\calQ}$;
we claim $\calY$ contains none of the maximal isotropic subspaces in the fibers of $\hat{q}$.  Indeed,
this is clear for the smooth fibers, which are disjoint from $\calY$.  The singular fibers
meet $\calY$ in {\em smooth} hyperplane sections.  However, given a smooth quadric hypersurface of even
dimension and a maximal isotropic subspace, any hyperplane section containing the subspace is necessarily
singular.  

This analysis implies that $\widetilde{\calW}_p$ has the geometric description given in the Lemma.  
Note that $\widetilde{\calW}\ra \calX$ is smooth, as each $\Lambda \subset \widetilde{\calQ}_p$ is a complete
intersection of (linear) Cartier divisors.  

Consider the mapping
$$\begin{array}{rcl}
\widetilde{\calW} & \ra & \widehat{\calW} \\
\Lambda & \mapsto & \gamma(\Lambda);
\end{array}
$$
we will show this is an isomorphism.   The key point is to demonstrate how the universal family of maximal isotropic
subspaces over $\widehat{\calW}$ induces a flat family of subschemes in $\widetilde{\calQ}$, corresponding to $\widetilde{\calW}$.
We apply a general observation:  Suppose we have a smooth variety $\calY$ embedded as a codimension-two
subvariety in a smooth variety $\calP$, and a flat family $\{\Lambda_t\}_{t\in T}$ of subvarieties
of $\calP$ over an integral base $T$ such that $\calY \cap \Lambda_t \subset \Lambda_t$ is either empty or Cartier for each $t\in T$.
Then the total transforms of the $\{\Lambda_t \}_{t\in T}$ in $\mathrm{Bl}_{\calY}\calP$ form a flat family
as well.  Indeed, the `valuative criterion for flatness' \cite[11.8.1]{EGAIV} reduces us to the
curve case, where dimensional considerations suffice \cite[III.9.7]{Hart}.  

The first modification yields morphisms
$$\begin{array}{rcl}
\widetilde{\calW}_i & \ra & \calW \\
\Lambda & \mapsto & \beta(\Lambda)
\end{array}
$$
of schemes smooth over $\Spec(\bZ)$.  Our description of the fibers $\widetilde{\calW}_p$ guarantees this
is bijective:  For each maximal isotropic subspace in a fiber of $q''$, there is a unique 
$\Lambda$ mapping to that subspace.  It follows that these morphisms are isomorphisms.

\end{proof}
The lemma implies the smoothness of $\calW \ra \calX$, which completes the proof of Proposition~\ref{proposition:Stein}.  
\end{proof}

\subsection{The variety of maximal isotropic subspaces of a quadric}
\label{s:maximal}

In this section, we explain how to compute explicit equations for the variety of one-dimensional isotropic (projective) subspaces over the universal quadric in $4$ variables. Let $\calH = \bP_\bZ^{9} = \Proj \bZ[a_{11},\dots,a_{44}]$, and let $\calQ \subseteq \bP^3_\bZ\times \calH$ be the universal family of quadrics, given by the vanishing of
\begin{equation}
\begin{split}
a_{11}X_1^2 &+a_{22}X_2^2+a_{33}X_3^2+a_{44}X_4^2 \\
&\qquad +a_{12}X_1X_2+a_{23}X_2X_3+a_{34}X_3X_4+a_{13}X_1X_3+a_{24}X_2X_4+a_{14}X_1X_4.
\end{split}
\label{eq:univ quadric}
\end{equation}

We define the degree-zero graded $\bZ$-module homomorphism
\[
\begin{array}{rcl}
\varphi\colon\bZ[X_1,X_2,X_3,X_4] &\to& \bZ[s,t] \\
X_i &\mapsto& p_is + p_i't \qquad i = 1,\dots,4.
\end{array}
\]
Consider $p_1,\dots,p_4'$ as indeterminates, subject to the restriction that the rank of the matrix
\[
\begin{pmatrix}
p_1 & p_2 & p_3 & p_4 \\
p_1' & p_2' & p_3' & p_4'
\end{pmatrix}
\]
is maximal. Then the variety of one-dimensional isotropic subspaces $\calL$ of $\calQ$ is given by the vanishing of $\calQ$ along the image of $\Sym^2\varphi$.  Concretely, substitute the expressions $X_i = p_is + p_i't$ for $i = 1,\dots,4$ in~\eqref{eq:univ quadric} and think of the result as the zero polynomial in the variables $s$ and $t$ with coefficients in $\bZ[a_{11},\dots,a_{44},p_1,\dots,p_4']$. Let $I$ be the ideal generated by these coefficients.

To give equations for $\calL$ as a subscheme of $\calH \times \bG(1,3)$, where we consider $\bG(1,3) \subseteq \Proj \bZ[p_{12},p_{13},p_{14},p_{23},p_{24},p_{34}]$ under the Pl\"ucker embedding, let
\[
\begin{split}
J &:=  \langle p_{12} - p_1p_2' + p_2p_1', p_{13} - p_1p_3' + p_3p_1', p_{14} - p_1p_4' + p_4p_1', \\
    &\qquad p_{23} - p_2p_3' + p_3p_2', p_{24} - p_2p_4' + p_4p_2', p_{34} - p_3p_4' + p_4p_3' \rangle
\end{split}
\]
be an ideal of the ring $\bZ[p_1,\dots,p_4',p_{12},\dots,p_{34}]$. Denote by $I^e$ and $J^e$ the extensions of $I$ and $J$, respectively, to the ring
\[
\bZ[p_1,\dots,p_4',a_{11},\dots,a_{44},p_{12},\dots,p_{34}].
\]
We obtain defining equations for $\calL$ by eliminating $p_1,\dots,p_4'$ from the ideal $I^e + J^e$.

Note that $\calL \times_{\calH} \calU=\calW$, the maximal isotropic subspaces of the quadrics with isolated singularities. To obtain equations for $\calW$, we focus on the distinguished affine subset of $\calU$ given by
\[
\frac{1}{2}\det
\begin{pmatrix}
2a_{11} & a_{12} & a_{13} \\
a_{12} & 2a_{22} & a_{23} \\
a_{13} & a_{23} & 2a_{33}
\end{pmatrix}=4a_{11}a_{22}a_{33}+a_{12}a_{13}a_{23}-(a_{11}a_{23}^2+a_{22}a_{13}^2+a_{33}a_{12}^2) \neq 0.
\] 
If $\Char k \neq 2$, then the result is as follows:
in the notation of Example~\ref{ex:discriminants}, let $\Delta = L^2 + 4M + 16N$, and write $K = k(\sqrt{\Delta})$ for the \'etale $k$-algebra $Z(Q)$. Letting
\[
\lambda(p_{14},p_{24},p_{34}) =  a_{11}p_{14}^2 + a_{12}p_{14}p_{24} + a_{13}p_{14}p_{34} + a_{22}p_{24}^2 + a_{23}p_{24}p_{34} + a_{33}p_{34}^2,
\]
the two components of $\calW$ over the field $K$ are given by 
\begin{align*}
\calW_1 &= \{\lambda = 0\} \cap \{p_{ij} = l_{ij}(p_{14},p_{24},p_{34}) : (ij) \in \{(12),(13),(23)\} \} \\
\calW_2 &= \{\lambda = 0\} \cap \{p_{ij} = l'_{ij}(p_{14},p_{24},p_{34}) : (ij) \in \{(12),(13),(23)\} \},
\end{align*}
where $l_{ij}$ and $l'_{ij}$ are linear forms with coefficients in $K$.

\section{From cubic fourfolds to K3 surfaces}
\label{S:CubicToK3}

In this section we explain how to explicitly construct, over an \emph{arbitrary} field, a (possibly singular) K3 surface
from a cubic fourfold containing a plane.

\subsection{Degree two K3 surfaces}
Let $(X,f)$ denote a polarized K3 surface $X$ of degree $2$ over a field $k$;
assume that $|f|$ is base-point free.  The 
induced morphism
\[
\phi_{|f|}\colon X \to \bP_k^2
\]
is a flat double cover and induces an exact sequence
$$0 \ra \calO_{\bP^2} \ra {\phi_{|f|}}_*\calO_X \ra Q \ra 0,$$
which is split by the trace mapping unless $\Char(k)=2$.  In any event,
$Q$ is torsion-free because $\calO_{\bP^2}$ is integrally
closed;  and it is reflexive of rank one, thus invertible.
As $H^2(\bP^2,Q) \simeq
H^2(X,\calO_X)\simeq k$, we have $Q\simeq \calO_{\bP^2}(-3)$ and
our exact sequence still splits, albeit noncanonically.

The surface $X$ then is embedded in $\bP(\calO_{\bP^2} \oplus \calO_{\bP^2}(-3))$
and disjoint from the distinguished section.  Blowing this down, we realize
$X$ as a smooth sextic hypersurface in the weighted projective space $\bP_k(1,1,1,3) = 
\Proj k[x,y,z,w]$:
\begin{equation}
\label{eq:K3}
w^2 + \alpha(x,y,z)w + \beta(x,y,z) = 0,
\end{equation}
where $\alpha,\beta \in k[x,y,z]$ have respective degrees $3$ and $6$. In these coordinates,
the map $\phi_{|f|}$ is the restriction to $X$ of the natural projection $\bP_{k}(1,1,1,3) \dasharrow \bP_{k}^2$.
If $\Char(k) \neq 2$, then completing a square we may (and do) assume that $\alpha = 0$. In 
this case, the ramification curve is $C = \{\beta(x,y,z) = 0\}$.
In characteristic $2$, \eqref{eq:K3} yields a singular surface
if $\alpha(x,y,z)=0$.  (Indeed, such a surface is an inseparable double cover
of $\bP^2$, typically with $21$ ordinary double points;  its
minimal resolution is a supersingular K3 surfaces \cite{Shimada}.)
The ramification locus of $\phi_{|f|}$ is a curve
$C = \{\alpha(x,y,z) = 0\}$.

\subsection{Cubic fourfolds containing a plane} 
From now on, assume that 
$(Y,P)$ is a smooth cubic fourfold containing a plane $P$, defined
over a field $k$. 
Projection from $P$ induces a quadric surface bundle
$q\colon \wY \ra \bP^2_k.$
Let $r\colon\calW \ra \bP^2_k$ denote the relative
variety of lines, parametrizing lines in $Y$ contained in fibers of $q$.  

We may assume, without loss of generality, that
\[
\Proj k[X_1,X_2,X_3]=P \subset Y \subset \bP^5_k=\Proj k[X_1,X_2,X_3,Y_1,Y_2,Y_3]
\]
so that $P=\{Y_1=Y_2=Y_3=0\}$. The defining equation for $Y$ takes the form
\begin{equation}
\label{eq:generalfourfold}
L_{11}X_1^2+\cdots + L_{33}X_3^2+Q_{14}X_1+Q_{24}X_2+Q_{34}X_3+C_{44}=0,
\end{equation}
where $L_{ij}$, $Q_{ij}$ and $C_{44}$ are homogeneous forms in $k[Y_1,Y_2,Y_3]$ of respective degrees $1$, $2$ and $3$.
The blow-up of $\bP^5_k$ along $P$ may be identified with the projective
bundle $\bP(E)$, where
$$E\simeq \calO_{\bP_k^2}^{\oplus 3} \oplus \calO_{\bP_k^2}(-1).$$
The quadric bundle $q\colon\wY \ra \bP_k^2$ factors through $\bP(E)$;
the defining equation for $\wY \subset \bP(E)$ is given by some
$s \in \Gamma(\Sym^2 E^{\vee})$. If $\Char(k) \neq 2$, then we
may interpret $s$ as a homomorphism
$s\colon E\ra E^{\vee}$ with $s^{\vee}=s$.  Using our trivialization for $E$, we can express
\begin{equation} \label{eqn:matrix}
s=\left( \begin{matrix} 2L_{11} & L_{12} & L_{13} & Q_{14} \\
		          L_{12} & 2L_{22} & L_{23} & Q_{24} \\
			  L_{13} & L_{23} & 2L_{33} & Q_{34} \\
			  Q_{14} & Q_{24} & Q_{34} & 2C_{44}
		\end{matrix} \right).
\end{equation}
In other words, $\wY$ is isomorphic to the subscheme of $\Proj k[X_1,X_2,X_3,X_4] \times \Proj k[x,y,z]$ given by
\begin{equation}
\label{eq:projection}
\begin{split}
L_{11}(x,y,z)X_1^2 &+ \cdots + L_{33}(x,y,z)X_3^2 + C_{44}(x,y,z)X^4  \\
&+ Q_{14}(x,y,z) X_1X_4 + Q_{24}(x,y,z)X_2X_4+Q_{34}(x,y,z)X_3X_4=0.
\end{split}
\end{equation}
and the bundle map $q\colon\wY \to \bP_k^2$ is the projection onto the second factor.
The discriminant curve in $\bP_k^2$ has equation $\{F=0\}$, where $F$ is the determinant of the matrix (\ref{eqn:matrix});
here $X$ satisfies the equation 
\begin{equation}
\label{eq:generalK3charnot2}
w^2-F=0
\end{equation}
in $\bP_k(1,1,1,3) = \Proj k[x,y,z,w]$.

On the other hand, if $\Char k = 2$, then~\eqref{eq:projection} still holds, though representation of $s$ by~\eqref{eqn:matrix} does not. In this case, Example~\ref{ex:discriminants} shows that
\begin{equation}
\label{eq:generalK3char2}
w^2 + L(x,y,z)w + M(x,y,z) = 0
\end{equation}
is an equation for $X \subset \bP_k(1,1,1,3)$, where
\begin{align*}
L &= L_{12}Q_{34}+L_{13}Q_{24}+L_{23}Q_{14}, \\
M &= -(L_{12}L_{23}Q_{34}Q_{14} + L_{13}L_{23}Q_{24}Q_{14}+L_{12}Q_{24}L_{13}Q_{34}) \\
    &\quad +(L_{11}L_{23}Q_{24}Q_{34}+L_{22}L_{13}Q_{34}Q_{14}+L_{33}L_{12}Q_{24}Q_{14}+C_{44}L_{12}L_{23}L_{13}) \\
    &\quad -(L_{11}L_{22}Q_{34}^2+L_{11}L_{33}Q_{24}^2+L_{11}C_{44}L_{23}^2+L_{22}L_{33}Q_{14}^2+L_{22}C_{44}L_{13}^2+L_{33}C_{44}L_{12}^2).
\end{align*}

We apply this construction in two examples that are instrumental in the proof Theorem~\ref{thm: Main}.

\begin{example}
\label{ex: Char 2}
The K3 surface $X\subset \bP_{\bF_2}(1,1,1,3)$ arising from the smooth cubic fourfold $\mathfrak{C}_2$ in $\bP^5_{\bF_2}$
\[
\begin{split}
X_1^2Y_2 &+ X_1^2Y_3 + X_1X_2Y_1 + X_1X_2Y_2 + X_1X_3Y_1 + X_1Y_2^2 + X_2^2Y_3 + X_2X_3Y_3 \\
&\quad + X_3^2Y_1 + X_3^2Y_3 + X_3Y_1Y_2 + X_3Y_2^2 + Y_1^2Y_3 + Y_1Y_3^2 + Y_3^3 = 0
\end{split}
\]
is given by
\begin{equation}
\label{eq: Char 2}
\begin{split}
w^2 &+ w(x^2y + y^3 + y^2z) \\ 
&\quad + x^5z + x^3y^2z + x^2y^3z + x^3yz^2 + x^2y^2z^2 + y^2z^4 + xz^5 + yz^5 + z^6 = 0.
\end{split}
\end{equation}
\end{example}

\begin{example}
\label{ex: Char 3}
 The K3 surface $X\subset \bP_{\bF_3}(1,1,1,3)$ arising from the smooth cubic fourfold $\mathfrak{C}_3$ in $\bP^5_{\bF_3}$
\[
2X_1^2Y_1 + X_1^2Y_3 + X_1X_3Y_2 + 2X_1Y_3^2 + X_2^2Y_3 + X_2Y_2^2 + X_3^2Y_1 + X_3Y_1^2 + 2X_3Y_2Y_2 + 2X_3Y_2^2 + 2Y_1^3
\]
is given by
\begin{equation}
\label{eq: Char 3}
\begin{split}
w^2 &= 2x^5z + x^4yz + x^4z^2 + 2x^3yz^2 + x^2y^4 + 2x^2y^3z + x^2y^2z^2 \\
		&\quad+ 2x^2yz^3 + xy^3z^2 + xy^2z^3 + 2xz^5 + y^6 + 2y^4z^2 + y^3z^3.
\end{split}
\end{equation}
\end{example}

\section{Constructing unramified Azumaya algebras}
 \label{S:constructing}

The Hodge-theoretic constructions of \S\ref{S:hodgetheory} can now be carried out using geometric techniques, valid over an arbitrary field.

\begin{theorem}
\label{theorem:getAzumaya}
Let $Y$ be a cubic fourfold smooth over a field $k$.  Suppose that
$Y$ contains a plane $P$, and let $\wY$ denote the blow up
of $Y$ along $P$, $q:\wY \ra \bP^2$ the corresponding quadric surface bundle, and
$r:\calW \ra \bP^2$ its relative variety of lines.
Assume that there exists no plane $P' \subset Y_{\bar{k}}$
such that $P'$ meets $P$ along a line.

Then the Stein factorization 
$$r: \calW \stackrel{\pi_1}{\ra} X \stackrel{\phi}{\ra} \bP^2$$
consists of a smooth $\bP^1$-bundle followed by a degree-two
cover of $\bP^2$, which 
is a K3 surface.
\end{theorem}

\begin{proof}
Our hypothesis means that no geometric fiber of $q$ contains a plane;
thus the fibers of $q$ have at worst isolated singularities.
All the statements, with the exception that $X$ is a K3 surface,
then follow directly from Proposition~\ref{proposition:Stein}. 
For the remaining claim, we may assume that $k$ is algebraically closed.

We start with an elementary result.

\begin{lemma} \label{lemma:singdisc}
Let $B$ be a smooth surface over an algebraically closed field of characteristic $\neq 2$.
Suppose that $q:Q\ra B$ is quadric surface bundle such that
\begin{itemize}
\item{the generic fiber is smooth;}
\item{the degenerate fibers have at worst isolated singularities.}
\end{itemize}
Then the discriminant curve $\Delta \subset B$ is smooth if and only if $Q$ is smooth.  
\end{lemma}
\begin{proof}
Let $b\in \Delta$ be a point in the discriminant curve;  we may
replace $B$ with a formal neighborhood of $b$.  Let $s$ and $t$ be local coordinates 
centered at $b$.  Set $Q_b=q^{-1}(b)$ with equation
$$x^2+y^2+z^2=0.$$
The local equation of $Q$ takes the form
$$x^2+y^2+z^2=a_{11}(s,t)+2a_{12}(s,t)x+2a_{13}(s,t)y+\cdots,$$
where the coefficients $a_{ij}(s,t)$ are power series with 
vanishing constant term.  This is singular at $x=y=z=s=t=0$ precisely
when the linear term of $a_{11}(s,t)$ vanishes.  
We have 
$$\det \left(
\left( \begin{matrix} 0 & 0 & 0 & 0 \\
				0 & 1 & 0 & 0 \\
                               0 & 0 & 1 & 0 \\
	  		       0 &  0 & 0 & 1
		\end{matrix} \right)+
\left( \begin{matrix} a_{11} & a_{12} & a_{13} & a_{14} \\
				a_{12} & a_{22} & a_{23} & a_{24} \\
                               a_{13} & a_{23} & a_{33} & a_{34} \\
	  		       a_{14} &  a_{24} & a_{34}  & a_{44} 
		\end{matrix} \right)\right)=a_{11}+\text{h.o.t.},
$$
which implies that discriminant vanishes to first order whenever $a_{11}$ vanishes
to first order.
\end{proof}

Lemma~\ref{lemma:singdisc} yields the following amplification of~\cite[Lemme2, p.\ 582]{Voisin}:
If $\Char k \neq 2$ then the discriminant curve $\Delta \subset \bP^2_k$ 
of $q$ is a smooth plane sextic curve.
Thus in this case $X$ is a K3 surface, because it is a double cover of $\bP^2_k$ ramified along a smooth sextic curve.

If $\Char k = 2$, then we argue as follows:  The surface $X$ can only be singular at points on the ramification locus of $\phi$. Choose local coordinates $s$ and $t$ at $b \in \bP^2_k$ in the branch locus of $\phi$. Set $\wY_b=q^{-1}(b)$ with equation
\[
xy + z^2=0.
\]
Then, the local equation of $\wY$ takes the form
\[
xy + z^2 + A(s,t) + B(s,t)z = 0,
\]
where $A(s,t)$ and $B(s,t)$ are power series with vanishing constant terms.  On the other hand, the local equation of $X$ is
\[
w^2 + B(s,t)w + A(s,t) = 0
\]
(see Example~\ref{ex:discriminants}). It is then easy to see that $X$ is singular at $s = t = w = 0$ precisely when $\wY$ is singular at $x = y = z = s = t = 0$. However, $\wY$ is smooth over $k$, and thus so is $X$. It follows that $X$ is K3, because it is isomorphic to a smooth sextic hypersurface in $\bP_k(1,1,1,3)$.
\end{proof}

\section{Computation of the Picard group}
\label{S:Picard}

\subsection{van Luijk's Method}
\label{ss: van Luijk}

To prove that the geometric Picard rank $\rho(X)$ of a K3 surface (over $\bQ$) arising from a cubic fourfold is one, we use a method due to van Luijk \cite{vanLuijk}---see~\cite{ElsenhansJahnelRefinement} for a recent refinement of this method. The idea is as follows: Let $(X,f)$ be a polarized K3 surface over a number field. Suppose that $X$ has good reduction modulo $p$ for two distinct primes $p_1$ and $p_2$, and assume that 
\begin{enumerate}
\item the Picard rank of each reduction is two---the rank must be even over a finite field; 
\item the discriminants of the Picard groups of the reductions are in different square class.
\end{enumerate}
Then $\rho(X) = 1$. To prove that condition~(2) holds, we construct \emph{explicit} full-rank sublattices of the respective Picard groups, whose discriminants are in different square classes. For condition~(1), we first obtain upper bounds for the Picard rank of a reduction as follows: Let $(\calX_p,\bar f)$ be a polarized K3 surface over a finite field $k$ of cardinality $q = p^r$ ($p$ prime), and let $\ell \neq p$ be a prime.  Let $\Phi\colon \calX_p \to \calX_p$ be the $r$-th power absolute Frobenius and write $\Phi^*$ for the automorphism of $H^2_{\et}(\calX_p,\bQ_\ell)$ induced by the action of $\Phi\times 1$ on $\calX_p\times \bar k$. Then the number of eigenvalues of $\Phi^*$ of the form $q\zeta$, with $\zeta$ a root of unity, gives an upper bound for $\rho(\calX_p)$~\cite[Corollary 2.3]{vanLuijk}.

We determine the eigenvalues of $\Phi^*$ of the required type by computing the characteristic polynomial $f_{\Phi^*}(t)$ of $\Phi^*$ via Newton's identity:
\[
f_{\Phi^*}(t) = t^{22} + c_1t^{21} + c_2t^{20} +\cdots + c_{22},
\]
where
\[
c_1 = -\Tr(\Phi^*)\qquad\text{and}\qquad -kc_k = \Tr(\Phi^{*k}) + \sum_{i=1}^{k-1}c_i\Tr(\Phi^{*k-i}).
\]
The traces of powers of the Frobenius maps are in turn calculated by tallying up the points of $\calX_p$ over finite extensions of the base field,  using the Lefschetz trace formula:
\[
\#\calX_p(\bF_{q^n}) = q^{2n} + \Tr(\Phi^{*n}) + 1.
\]
The Weil conjectures give the functional equation
\[
p^{22}f_{\Phi^*}(t) = \pm t^{22}f_{\Phi^*}(p^2/t).
\]
Hence, to compute $f_{\Phi^*}(t)$, it may suffice to compute $\#\calX_p(\bF_{q^n})$ up to some $n < 22$. Typically, $n = 11$ is sufficient (if $c_{11} \neq 0$, giving a positive sign for the functional equation).  Sometimes, however, it is possible to do better, and get away with $n = 9$: see~\cite[\S 4]{ElsenhansJahnelrankone}.

\subsection{Special K3 surfaces of degree $2$}

\subsubsection{Characteristic $2$}
\label{ss: Char 2}

Let $(X,f)$ be a polarized K3 surface
of degree $2$ over an algebraically closed field $k$ with $\Char(k) = 2$, given in the form~\eqref{eq:K3}. Note that $\alpha \neq 0$ since $X$ is
smooth. The rational map
\[
\begin{array}{rcl}
\bP_k(1,1,1,3) &\dasharrow& \bP_k(1,1,1,3) \\
 {[x:y:z:w]} &\mapsto& [x:y:z:w + \alpha]
\end{array}
\]
restricts to an automorphism $\psi$ of $X$. 

Suppose that $X$ contains a divisor of the form
\begin{equation}
\label{eq:char2divisor}
\Gamma := \{l(x,y,z) = w + c(x,y,z) = 0\},
\end{equation}
where $l, c \in k[x,y,z]$ have respective degrees $1$ and $3$. Let $\Gamma'$ be
the scheme-theoretic image of $\Gamma$ under $\psi$. The intersection number $(\Gamma,\Gamma')$ is equal
to
\[
\deg(\Gamma\cap\Gamma') = \deg\big( \Proj k[x,y,z,w]/(l,\alpha,w + c)\big) =
\deg \big(\Proj k[x,y,z]/(l,\alpha)\big) = 3.
\] 
Since $\{ x = 0\} \in |f|$, it follows by similar computations that
\[
(\Gamma,f) = (\Gamma',f) = 1,
\]
and hence $\Gamma$ and $\Gamma'$ are irreducible, since $f$ is ample.
Let $V := \{l(x,y,z) = 0\} \subseteq \bP_{k}^2$. Then, on the one hand, we have
 $\phi^*V \in |f|$, so that $(\phi^*V,f) = 2$. On the other hand 
 $\phi^*V - \Gamma - \Gamma'$ is an effective divisor, and
 \[
 (\phi^*V - \Gamma - \Gamma',f) = 0;
 \]
 whence $\phi^*V = \Gamma + \Gamma'$, because $f$ is ample. Since $\psi$
 preserves the intersection form on $X$, the equality $(\phi^*V,\phi^*V) = 2$
 now shows that $(\Gamma,\Gamma) = (\Gamma',\Gamma') = -2$. We conclude that the classes
 of $\Gamma$ and $\Gamma'$ are independent in $\Pic X$ since the determinant 
 of the intersection matrix
 \[
 \bigg|
 \begin{matrix}
 -2 & \ 3 \\ 
 \ 3 & -2
 \end{matrix}
 \bigg| = -5
 \]
 is nonzero.  In particular, $\rho(X) \geq 2$.
 
 \begin{example}
 \label{ex:char2}
 Consider the K3 surface $X$ of Example~\ref{ex: Char 2}, defined by~\eqref{eq: Char 2}. The curve
\[
\Gamma := \{w = z = 0\}
\]
lies on this surface, and hence so does the image $\Gamma'$ of $\Gamma$ under $\psi$:
\[
\Gamma' = \{w +  x^2y + y^3 + y^2z = z = 0\}
\]
The classes of these curves over $\overline{\bF}_2$ together span a sublattice of $\Pic\overline{X}$ of discriminant $-5$, and thus $\rho(X) \geq 2$.
 \end{example}
  
 \subsubsection{Characteristic $\neq 2$}
\label{ss: Char not 2}

Let $(X,f)$ be a polarized K3 surface of degree $2$ over an algebraically closed field of $k$ with $\Char(k) \neq 2$, given in the form~\eqref{eq:K3} with $\alpha = 0$. 

Suppose that there exists a smooth conic $C'\subseteq \bP^2_k$ tangent at six (not necessarily distinct) points to the ramification curve $C$ of $\phi$. Arguing as in~\cite[\S 2]{ElsenhansJahnelrankone}, we have $\phi^*C = C_1 + C_2$, where $C_i$ is an irreducible divisor for $i = 1,2$, and the intersection matrix of the lattice spanned by the classes of $C_1$ and $C_2$ in $\Pic X$ has determinant
 \[
 \bigg|
 \begin{matrix}
 -2 & \ 6 \\ 
 \ 6 & -2
 \end{matrix}
 \bigg| = -32 \neq 0.
 \]
Hence $\rho(X) \geq 2$ in this case.

\begin{example}
\label{ex:char3}
The polynomial
\begin{equation}
\label{eq: C in Char 3}
\begin{split}
f(x,y,z) &= 2x^5z + x^4yz + x^4z^2 + 2x^3yz^2 + x^2y^4 + 2x^2y^3z + x^2y^2z^2 \\
		&\quad+ 2x^2yz^3 + xy^3z^2 + xy^2z^3 + 2xz^5 + y^6 + 2y^4z^2 + y^3z^3
\end{split}
\end{equation}
defines a smooth sextic curve $C \subseteq \bP_{\bF_3}^2$; the double cover $X$ of $\bP^2_{\bF_3}$ ramified along $C$ is the K3 surface $X$ of Example~\ref{ex: Char 3}, given by~\eqref{eq: Char 3}. The smooth conic $C'$ in $\bP_{\bF_3}^2$ given by
\[
2x^2 + 2xy + xz + 2y^2 = 0.
\]
is tangent to $C$ at six points: indeed, let $u$ be a multiplicative generator for $\bF_{3^2}$. Then $C'$ is parametrized by the map $\bP^1_{\bF_{3^2}} =\Proj(\bF_{3^2}[s,t]) \to \bP^2_{\bF_{3^2}}$ given by
\[
x = u^2t^2,\quad
y = u^6st,\quad\text{and}\quad
z = u^2s^2 + u^6st + u^2t^2.
\]
Substituting these expressions into $f(x,y,z)$ we obtain
\[
t^2(s^5 + s^4t + s^3t^2 + s^2t^3 + 2st^4 + t^5)^2
\]
The zeroes of this polynomial over $\overline{\bF}_3$ correspond to the geometric points of $C\cap C'$; they clearly all have multiplicity greater than one. We conclude that $\rho(X) \geq 2$.
\end{example}

\subsection{A K3 surface of degree $2$ over $\bQ$ with geometric Picard rank one}
We collect the results of this section in a single proposition for future reference.  The proposition is a simple application of van Luijk's method.

\begin{proposition}
\label{P: rank 1}
Let $X$ be a K3 surface of degree $2$ over $\bQ$. Suppose that the reductions $\calX_2$ and $\calX_3$ of $X$ at $p = 2$ and $3$ are isomorphic to the K3 surfaces~\eqref{eq: Char 2} and~\eqref{eq: Char 3}, respectively. Then $\rho(X) = 1$.
\end{proposition}

\begin{proof}
We use the notation introduced in \S\ref{ss: van Luijk}.  Counting $\bF_{p^n}$ points on $\calX_p$ for $n = 1,\dots,12$ and $p = 2, 3$, and using Newton's identities, we may compute the first twelve coefficients $c_1,\dots,c_{12}$ of the characteristic polynomial $f_{\Phi^*}(t)$ for $p = 2$ and $3$. We obtain:
$$
\begin{array}{c||c|c|c|c|c|c|c|c|c|c|c|c}
p & c_1 & c_2 & c_3 & c_4 & c_5 & c_6 & c_7 & c_8 & c_9 & c_{10} & c_{11} & c_{12} \cr
\hline
2 & -2 & -2 & 4 & 8 & -32 & 0 & 64 & 128 & -512 & 512 & 0 & 2\,048 \cr
\hline
3 & -1 & -18 & 9 & 135 & 162 & -243 & -3\,645 & -6\,561 & 26\,244 & 118\,098 & 0 & -106\,288
\end{array}
$$
(see Remark~\ref{rem:counts} below for the actual point counts).  Note that to determine the sign of the functional equation of the characteristic polynomial, we must compute $\#\calX_p(\bF_{p^{12}})$, since $c_{11} = 0$ in both cases. We now have enough information to determine the remaining coefficients of the characteristic polynomial.
Let $\tilde{f}_{p}(t) = p^{-22}f_{\Phi^*}(pt)$, so that the number of roots of $\tilde{f}_{p}(t)$  that are roots of unity gives an upper bound for $\rho(\calX_p)$. We have the following factorizations into irreducible factors
\begin{align*}
\tilde{f}_{2}(t) &= \frac{1}{2}(t-1)^2 (2t^{20} + 2t^{19} + t^{18} + t^{17} + 2t^{16} + t^{15} + t^{12} + t^8 + 
        t^5 + 2t^4 + t^3 + t^2 + 2t + 2), \\
\tilde{f}_{3}(t) &= \frac{1}{3}(t-1)(t+1) (3t^{20} - t^{19} - 3t^{18} + 2t^{16} + 2t^{15} + t^{14} - 3t^{13} - 2t^{12} \\
&\hspace{1.6in} + t^{11} + 4t^{10} + t^9 - 2t^8 - 3t^7 + t^6 + 2t^5 + 2t^4 - 3t^2 - t + 3).
\end{align*}
In both cases, the roots of the degree $20$ factor of $\tilde{f}_{p}(t)$ are not integral, so they are not roots of unity.  We conclude that $\rho(\calX_p) \leq 2$ for $p = 2, 3$, and hence that $\rho(X_p) = 2$, by our work in Examples~\ref{ex:char2} and~\ref{ex:char3}.

Finally, the discriminants of the full sublattices of $\Pic(\calX_p)$ exhibited in Examples~\ref{ex:char2} and~\ref{ex:char3} are in different square classes, and hence $\rho(X) = 1$.
\end{proof}

\begin{remark}
\label{rem:counts}
We record below the number of $\bF_{p^n}$-points of $\calX_p$. Set $N_n := \#\calX_p(\bF_{p^n})$; the second row of the table shows the counts for $p = 2$; the last row shows the counts for $p = 3$.
{\tiny
$$
\begin{array}{r|r|r|r|r|r|r|r|r|r|r|r}
N_1 & N_2 & N_3 & N_4 & N_5 & N_6 & N_7 & N_8 & N_9 & N_{10} & N_{11} & N_{12} \cr
\hline
7 & 25 & 73 & 249 & 1137 & 4273 & 16737 & 65313 & 264385 & 1047745 & 4203393 & 16767105 \cr
\hline
11&119&758 & 6707 & 58421 & 529472 & 4784357 & 43059323& 387449246 & 3486568169 & 31380849731 & 282429079832
\end{array}
$$}
Counting points over $\bF_{3^{11}}$ and $\bF_{3^{12}}$ is computationally time consuming. See \S\ref{Section: computations} for details.
\end{remark}

\section{Proof of Theorem~\ref{thm: Main}}
\label{S: proof}

\subsection{Geometric Picard rank one}
The K3 surface $X$ of Theorem~\ref{thm: Main} arises from the cubic fourfold $Y$ given by~\eqref{eq:fourfold}: this is a direct computation applying the construction of \S\ref{S:CubicToK3}.

\begin{remark}
An \emph{a posteriori} analysis shows that $Y$ does not contain a plane incident to $P$ along a line.  If it were to, then the quadric surface bundle $q\colon\wY \to \bP^2_{\bQ}$ would have at least one singular geometric fiber containing a plane, and an analysis analogous to that in the proof of Lemma~\ref{lemma:singdisc} would then show that the discriminant curve in $\bP^2_{\bQ}$ is singular. However, it is easy to check that this is not the case.
\end{remark}

\begin{proposition}
\label{prop:picardOne}
We have $\rho(X) = 1$.
\end{proposition}
\begin{proof}
The reductions of the cubic fourfold $Y$ at $p = 2$ and $3$ coincide with the cubic fourfolds $\mathfrak{C}_2$ and $\mathfrak{C}_3$ of Examples~\ref{ex: Char 2} and~\ref{ex: Char 3}, respectively. Thus, $X$ has good reduction at $p = 2$ and $3$, and its reductions at these primes are isomorphic to the K3 surfaces~\eqref{eq: Char 2} and~\eqref{eq: Char 3}, respectively. The claim now follows directly from Proposition~\ref{P: rank 1}.
\end{proof}

\subsection{The quaternion algebra}
\label{ss:quaternion}

Let $r\colon\calW \ra \bP^2_\bQ$ denote the relative
variety of lines, parametrizing lines in $Y$ contained in fibers of $q\colon\wY \to \bP^2_\bQ$. By Theorem~\ref{theorem:getAzumaya}, the Stein factorization of $r$ is given by 
\[
\calW \xrightarrow{\pi_1} X \xrightarrow{\phi} \bP^2_\bQ,
\]
where $\pi_1$ is a smooth $\bP^1$-bundle over $X$. The generic fiber of $\pi_1$ is a Severi-Brauer variety $\calC$ over the function field $k(X)$ of $X$. By our discussion in \S\ref{s:maximal}, $\calC$ is explicitly given by the conic
\[
(2x + 3y + z)p_{14}^2 + (3x + 3y)p_{14}p_{24} + (3x + 4y)p_{14}p_{34} + zp_{24}^2 + 3zp_{24}p_{34} + (x+3z)p_{34}^2 = 0
\]
 in $\Proj k(X)[p_{14},p_{24},p_{34}]$.  Completing squares and renormalizing we obtain
\[
p_{14}^2 = \alpha p_{24}^2 + \beta p_{34}^2,
\]
where
\begin{align*}
\alpha &= \frac{9x^2 + 18xy - 8xz + 9y^2 - 12yz - 4z^2}{4(2x + 3y + z)^2},\ \text{and} \\
\beta &= -\frac{9x^3 + 18x^2y + x^2z + 9xy^2 + 3xyz - 10xz^2 + 7y^2z - 9yz^2 - 3z^3}{(2x + 3y + z)(9x^2 + 18xy - 8xz + 9y^2 - 12yz - 4z^2)}.
\end{align*}
The conic $\calC$ is associated to the quaternion algebra $(\alpha,\beta) \in \Br(k(X))[2]$~\cite[Corollary~5.4.8]{GilleSzamuely}. Theorem~\ref{theorem:getAzumaya} can be recast as saying that $(\alpha,\beta)$ is in the image of the natural injection 
\[
\Br(X) \hookrightarrow \Br(k(X)).
\]

\subsection{Evaluation of the local invariants}

Since the K3 surface $X$ is projective, the natural inclusion $X(\bA_\bQ)\subseteq \prod_{p\leq \infty}X(\bQ_p)$ is a bijection. Class field theory gives rise to the constraint
\[
X(\bQ) \subseteq X(\bA_\bQ)^{\Br} := \bigg\{ (x_p)_p \in X(\bA_\bQ) \ {\big|}\ \sum_{p\leq \infty}\inv_p\big(\mathscr{A}(x_p)\big) = 0\text{ for every }\mathscr{A}\in \Br(X) \bigg\},
\]
where $\mathscr{A}(x_p) := \mathscr{A}_{x_p}\otimes_{\mathscr{O}_{X,x_p}}\bQ_p$ and $\inv_p\colon \Br\bQ_p \to \bQ/\bZ$ is the local invariant map for each $p \leq \infty$. 


A naive search for points on $X$ reveals the rational point $P_1 := [15:15:16:13\,752]$ and the real point $P_2 := [1:0:1:\sqrt{8}]$. Let $\mathscr{A} = (\alpha,\beta)$ be the quaternion algebra above, considered as an element of $\Br(X)$. Then
\[
\mathscr{A}(P_1) = \bigg({\frac{2\,276}{4\cdot 91^2}},-\frac{75\,852}{91\cdot 2\,276}\bigg)\quad\text{and}\quad 
\mathscr{A}(P_2) = \bigg({-\frac{3}{36}},-\frac{1}{3}\bigg)
\]
and thus
\begin{equation}
\label{eq:invariants}
\inv_\infty\big(\mathscr{A}(P_1)\big) = 0 \quad\text{and}\quad 
\inv_\infty\big(\mathscr{A}(P_2)\big) = \frac{1}{2}.
\end{equation}

\begin{proof}[Proof of Theorem~\ref{thm: Main}]
The K3 surface $X$ has geometric Picard rank one by Proposition~\ref{prop:picardOne}. Thus, the algebraic Brauer group of $X$ is trivial, because it is isomorphic, by the Hochschild-Serre spectral sequence, to the Galois cohomology group
\[
H^1(\Gal(\overline{\bQ}/\bQ),\Pic(\overline{X})),
\]
and the Galois action is trivial on the rank-one free $\bZ$-module $\Pic(\overline{X})$.
 Therefore, the quaternion algebra $\calA = (\alpha,\beta)$ of \S\ref{ss:quaternion}, is a transcendental element of $\Br(X)$.

Define $\calP \in X(\bA_\bQ)$ as the point that is equal to $P_1$ at all finite places, and is $P_2$ at the real place. By~\eqref{eq:invariants}, it follows that
\[
\sum_{p\leq \infty}\inv_p\big(\mathscr{A}(\calP_p)\big) = \frac{1}{2},
\]
and thus $\calP \notin X(\bA_\bQ)^{\Br}$. This shows that $X$ does not satisfy weak approximation.
\end{proof}

\begin{remarks}
\ 
\begin{enumerate}
\item It is natural to ask if the methods of this paper can also yield counter-examples to the Hasse principle.  Numerical experiments suggests that, in the absence more constraints, the primes of bad reduction for the $K3$ surfaces we construct give rise to only mild singularities.  It is possible that such singularities cannot force the evaluation maps at these bad places to be constant. Thus, more constraints on the initial cubic fourfold may be necessary to violate the Hasse principle with a single Brauer element. 

We also remark that imposing conditions on a smooth cubic fourfold at two different primes tends to create rather large primes of bad reduction for the associated K3 surface.  Calculating the evaluation maps at these primes can be a computational challenge. The methods of~\cite{ElsenhansJahnelOnePrime} may well obviate this difficulty.

\item It would be interesting to figure out if $\overline{X(\bQ)} = X(\bQ)^{\Br}$. Such an equality holds for many del Pezzo surface. However, in our case, we do not even know if the set $X(\bQ)$ is Zariski dense!
\end{enumerate}
\end{remarks}

\section{Computations}
\label{Section: computations}

In this section we briefly outline the steps used to construct the example in Theorem~\ref{thm: Main}, for the benefit of those wishing to carry out similar computations. \\

\noindent{\bf Step 1:} Generate a list of homogeneous random polynomials in $\bF_p[x,y,z]$ ($p = 2$ or $3$)
\[
\{L_{11},L_{12},L_{13},L_{22},L_{23},L_{33},Q_{14},Q_{24},Q_{34},C_{44}\}
\]
where the $L_{ij}$ are linear, the $Q_{ij}$ are quadratic and $C_{44}$ is cubic. Check smoothness of the cubic fourfold~\eqref{eq:generalfourfold} and the K3 surface $\calX_p$ given by~\eqref{eq:generalK3charnot2} (or~\eqref{eq:generalK3char2} if $p = 2$). Start over if either one is not smooth.

\medskip

\noindent{\bf Step 2:} If $p = 2$, then test to see if the K3 surface has a divisor of the form~\eqref{eq:char2divisor}; if $p = 3$, then use \cite[Algorithm~12]{ElsenhansJahnelrankone} to test for the existence of a conic tangent at six points to the branch curve of $\calX_p$. In each case, go back to Step 1 if the test fails to reveal curves of the required type.

\medskip

\noindent{\bf Step 3:} If $p = 2$, then count $\bF_{2^n}$-points on the K3 surface for $n = 1,\dots, 12$. Use this to compute the first $12$ coefficients $c_1\dots,c_{12}$ of the characteristic polynomial of the action of Frobenius on $H^2_{\et}(\calX_p,\bQ_\ell)$, as described in~\S\ref{ss: van Luijk}. If either $c_{11} \neq 0$ or $c_{12}\neq 0$ then compute the remaining coefficients using the functional equation for the characteristic polynomial. Use this information to obtain an upper bound $\rho_{\text{up}}$ for the geometric Picard rank of the K3 surface. If $\rho_{\text{up}} > 2$, then go back to Step 1.

If $p = 3$, then use~\cite[Algorithm 15]{ElsenhansJahnelrankone} to count $\bF_{3^n}$-points on the K3 surface for $n = 1,\dots, 10$. This algorithm uses the fact that $\calX_p$ is defined over $\bF_p$, and counts \emph{Galois orbits} of points. This saves a factor of $n$ when counting $\bF_{p^n}$-points. Next, use \cite[Algorithms~21 and~23]{ElsenhansJahnelrankone} to determine (with a very high degree of certainty) an upper bound $\rho_{\text{up}}$ for the geometric Picard rank of the K3 surface. If $\rho_{\text{up}} > 2$ then go back to Step 1. If $\rho_{\text{up}} = 2$ then count $\bF_{3^n}$-points on the K3 surface for $n$ large enough to ascertain the characteristic polynomial of Frobenius with total certainty; this last step may not be necessary.

\medskip

\noindent{\bf Step 4} Use the Chinese remainder theorem to construct a cubic fourfold over $\bQ$ that reduces, mod $2$ and $3$, to the cubic fourfolds of Step 1. The resulting K3 surface $X$ reduces to the corresponding K3 surfaces over $\bF_2$ and $\bF_3$, and thus $\rho(X) = 1$. Compute the quaternion algebra $\mathscr{A} \in \Br(X)$, and perform a naive search for rational and real points to test for an obstruction to weak approximation, as in the Proof of Theorem~\ref{thm: Main}. \\

For Step 3, we implemented~\cite[Algorithm~15]{ElsenhansJahnelrankone} in {\tt Magma} to count points on double covers of the projective plane over $\bF_{3^n}$ for $n = 1,\dots,10$.  To count points over $\bF_{3^{11}}$ and $\bF_{3^{12}}$ in a reasonable amount of time, we programmed the arithmetic of these fields in {\tt C++}. We obtained significant speed gains by implementing the following ideas:

\begin{itemize}
\item We choose an isomorphism $\bF_{3^n} \cong \bF_3[X]/(f)$, where $f$ is an irreducible polynomial of degree $n$ with a minimal number of nonzero coefficients, and store elements of $\bF_{3^n}$ using representatives of degree less than $n$. 

\item For each element~$g \in \bF_{3^n}$, we precalculate values of $g^k$, $2 \leq k \leq 6$ (use H\"orner's algorithm to perform exponentiation). Intermediate results give polynomials in $X$ of degree less than $2n - 1$; we reduce these polynomials modulo $f$ inline, rather than using pre-stored tables. As a by-product, we build a list of squares in $\bF_{3^n}$, which allows us to quickly answer the query ``is $g$ a square?''.

\item To evaluate a polynomial such as~\eqref{eq: C in Char 3}, on the affine patch $x = 1$, we compute each monomial by multiplying two elements from our precalculated list, obtaining a polynomial in $X$ of degree less than $2n-1$. We add these polynomials in $X$ \emph{first}, and then perform a single reduction mod $f$. If the result, for example, is a nonzero square, then we have two $\bF_{3^n}$-rational points on the corresponding K3 surface~\eqref{eq: Char 3}.
\end{itemize}

%




\bibliographystyle{alpha}

\bibliography{cubicK3}

\end{document}